\documentclass[12pt]{article}

\usepackage{amsfonts, amssymb, amsmath, amsthm, epsf, graphicx, cite}

\usepackage{fullpage}

\newcommand{\cP}{{\mathcal P}}
\newcommand{\cR}{{\mathcal R}}

\newcommand{\cW}{{\mathcal W}}
\newcommand{\ux}{{\underline  x}}
\newcommand{\ox}{{\overline x}}

\newcommand{\uU}{{\underline U}}

\newcommand{\oF}{{\overline F}}

\newcommand{\bbN}{{\mathbb N}}
\newcommand{\bbR}{{\mathbb R}}

\newtheorem{theorem}{Theorem}[section]
\newtheorem{lemma}{Lemma}[section]

\theoremstyle{definition}
\newtheorem{definition}{Definition}[section]
\newtheorem{remark}{Remark}[section]

\numberwithin{equation}{section} \numberwithin{figure}{section}

\title{Asymptotics of sign-changing patterns in hysteretic systems with diffusive thresholds}

\author{Pavel Gurevich\footnote{Free University of Berlin, Institute for Mathematics,
Arnimallee 3, 14195 Berlin, Germany \& Peoples' Friendship University of Russia, Miklukho-Maklaya str.~6, 117198 Moscow, Russia;  email: gurevich@math.fu-berlin.de},\qquad
Dmitrii Rachinskii\footnote{Department of Mathematical Sciences,
The University of Texas at Dallas,
800 West Campbell Road,
Richardson, TX 75080-3021, USA; tel. +1-972-883-4401;\ fax: +1-972-883-6622;\ email:
dmitry.rachinskiy@utdallas.edu\ (corresponding author)}}

\begin{document}

\maketitle


\begin{abstract}
We consider a reaction-diffusion system including discontinuous hysteretic relay operators in reaction terms.
This system is motivated by an epigenetic model that describes the evolution of a population of organisms
which can switch their phenotype in response to changes of the state of the environment.
The model exhibits formation of patterns in the space of distributions of the phenotypes
over the range of admissible switching strategies. We propose asymptotic formulas for the pattern and the process
of its formation.
\end{abstract}

\bigskip


\section{Introduction}
In this paper, we consider a reaction-diffusion equation that includes discontinuous two-state two-threshold relay operators.
This equation was proposed as a model for dynamics of a population of bacteria that can switch their phenotype in response to
changes of the environment \cite{gena1, GurRachDiffThresh2013}. The model exhibits formation of patterns
in the space of distributions of the two phenotypes over an admissible range of the switching thresholds.
The convergence of solutions to a stationary pattern observed numerically in \cite{gena1} was proved in \cite{GurRachPattern2015}.
The objective of this work is to obtain asymptotic formulas for the pattern and the timing of its formation
under the assumption of slow diffusion.

The model includes several time scales. The ``spatial'' variable is a scalar parameter $x$ that measures the distance between the switching thresholds
of a bacterium  and characterizes
its switching strategy (the smaller the $x$ the more responsive is the strategy to variations of the environmental conditions).
That is, a solution describes the evolution of the distribution of each of the two phenotypes over possible switching strategies
parameterized by $x$. The diffusion of this distribution models sporadic changes of the switching strategy by individual organisms; it is assumed to be slow
compared to the rate of growth and competition precesses. The rate of transition from one phenotype to the other in bacteria
is assumed to be much higher than the rate of other processes. Every such transition is modeled by an instantaneous switch
of a relay from state $1$ to state $-1$ or vice versa (where each state represents a particular phenotype).

The central assumption that we make is that the switching threshold for the transition of a relay (bacterium) from state $1$
to state $-1$ is different from the switching threshold for the opposite transition (the difference of the thresholds $2x$ is called the bi-stability range).
That is, we assume hysteresis, or multi-stability, in the switching response of bacteria to the exogenous stimulus.
This hysteresis
can be viewed as form of the persistent memory of past environmental conditions by organisms.
Since the seminal work of  Max Delbr\"uk \cite{delb},
epigenetic differences arising in the process of
cell differentiation have been attributed to multi-stability of living forms. A
classical example of such multi-stability is the behavior of {\em lac-operon} in {\em E.~coli}
({\em lac-operon} is a collection of genes associated with transport
and metabolism of lactose in the bacterium; expression of these genes can
be turned on by molecules called inducers).
Experimental studies of the regulation of enzymes in {\em E.~coli} and yeast that were conducted as early as in the 1950-60's
 effectively demonstrated that two phenotypes each
associated with ``on'' and ``off'' state of {\em lac-operon}
expression can be obtained from the same culture. The fraction of
each phenotype in the total population depended on the history of
exposure to the inducer, and  each phenotype remained stable
through multiple generations of the bacteria \cite{benz, cohn,
cohn1, mono, novi, spie, wing}. This behavior resembles the
two-threshold hysteretic relay shown in Fig.~\ref{hyst}. Later,
multi-stable gene expression and hysteresis have been
well-documented in many natural as well as artificially
constructed systems \cite{3, 15, 17, 5, 10}. In particular, recent
experiments using molecular biology methods permitted a further
study of the region of bi-stability of the {\em lac-operon} when
multiple input variables are used to switch the {\em lac-operon}
genes on and off.

A substantial amount of experimental and theoretical work provides
an evidence that organisms use various (often randomized)
strategies of phenotype switching in order to adapt to changing
environmental conditions (see  \cite{KL} and references therein).
It is intuitively clear that phenotypic diversity within the
population can help to increase the chances of survival in varying
environment. Switching strategies that establish phenotypic
diversity are called bet-hedging \cite{1}.
Indeed, experimental and theoretical models of adaptation demonstrate that bet-hedging can evolve to maximize the net growth of the total population \cite{to1, to2, to3, to4, to5, to6}.
In \cite{gena}, a simple differential model  was used to show that switching strategies with two thresholds
where the switching moment depends not only on the state of the environment, but also on the phenotype itself,
can
further increase Darwinian fitness of species. These strategies
exploit hysteresis and are described by a two-threshold relay
shown in Fig.~\ref{hyst} in the adiabatic limit when switching is
faster than the characteristic rate of the environment variations.
It was also shown in \cite{gena} that the optimal bi-stability
range, that is the optimal separation of thresholds $2x$ of the
relay, which maximizes the net growth rate of the population is
different for different environmental inputs. For this reason, the
model that is considered in this work includes the distribution of
bacteria (relays) over a set of switching strategies which have
different separation of switching thresholds with $x$ taking the
values from an interval $[\underline{x},\overline{x}]$. The
diffusion process included in the model acts as a factor that
diversifies the switching thresholds (strategies) in the
population, while the competition provides a mechanism of
selection that may favor a subpopulation with a specific $x$ for a
given law of variation of the state of environment. More detailed
discussion of the model assumptions and the biological background
can be found in \cite{gena1, GurRachDiffThresh2013, gena}.

Parabolic equations with distributed non-ideal relays were
previously used for modeling spatial structures that were observed
in spatially distributed colonies of bacteria, for example concentric rings in
Petri dish experiments \cite{Jaeger, Jaeger1, Jaeger2, Gurevich, Gurevich1, Japan}. In
these models, all the relays have the same fixed separation of
thresholds. The mechanism of pattern formation discussed in this
work is based on the long-term memory of the hysteretic relays and is related to dynamics of free boundary
separating the domain where relays are in state $1$ from the domain where relays are in state $-1$
on the interval $[\underline{x},\overline{x}]$.
This mechanism
is different from the mechanism found in \cite{Jaeger, Jaeger1,
Jaeger2}. Other important aspects of dynamics of parabolic and
hyperbolic differential equations with distributed relays were
studied in \cite{BrokateBook, Botkin, term0, term1, term2, colli,
Kopfova, Kopfova1, KrejciBook,  Visintin, VisintinBook} in
relation to multiple applied problems (not related to patterns)
such as phase transitions, population dynamics,  multi-phase flows
in porous media, thermostat control and others.

In the following Sections 2 and 3 we present the model and a statement on its well-posedness that has been obtained earlier.
Section 4 contains the main statement on the asymptotics of the patterns and its proof.

\section{Model description}\label{secModelDescription}

In this paper, we consider a class of models, which attempt to
account for a number of phenomena, namely (a) switching of
bacteria between two phenotypes (states) in response to variations
of environmental conditions; (b) hysteretic switching strategy
(switching rules) associated with bi-stability of phenotype
states; (c) heterogeneity of the population in the form of a
distribution of switching thresholds; (d) bet-hedging in the form
of diffusion between subpopulations characterized by different
bi-stability ranges; and, (e) competition for nutrients. The
resulting model is a reaction-diffusion system including, as
reaction terms, discontinuous hysteresis relay operators and the
integral of those. This integral can be interpreted as the
Preisach operator~\cite{VisintinBook, pr1, pr2, pr3, pr4, pr5, pr6, pr7, pr8, pr9} with a time dependent
density (the density is a component of the solution describing the
varying distribution of bacteria). In~\cite{GurRachPattern2015},
we have shown that fitness, competition and diffusion can act
together to select a nontrivial distribution of phenotypes
(states) over the population of thresholds. The main goal of this
paper is to give a quantative description of this distribution for
small diffusion.

We assume that each of the two phenotypes, denoted by $1$ and
$-1$, consumes a different type of nutrient (for example, one
consumes lactose and the other glucose). The amount of nutrient
available for phenotype $i$ at the moment $t$ is denoted by
$f_i(t)$ where $i=\pm 1$. The model is based on the following
assumptions (see~\cite{GurRachDiffThresh2013} for further
discussion).

\begin{enumerate}
\item Each bacterium changes phenotype in response to the
variations of the variable $w=f_1/(f_1+f_{-1})-1/2$, which
measures the deviation of the relative concentration of the first
nutrient from the value 1/2 in the mixture of the two nutrients.

\item The input $w=w(t)$ is mapped to the (binary)   phenotype
(state) of a bacterium $r(t)={\mathcal R^x}(w)(t)$, where
${\mathcal R^x}$ is the non-ideal relay operator   with symmetric
switching thresholds $x$, $y=-x$ with $x>0$; see Fig.~\ref{hyst}
and the rigorous definition~\eqref{nonideal} in
Section~\ref{secSetting}.
\begin{figure}[ht]
  \centering
  \includegraphics[width=90mm]{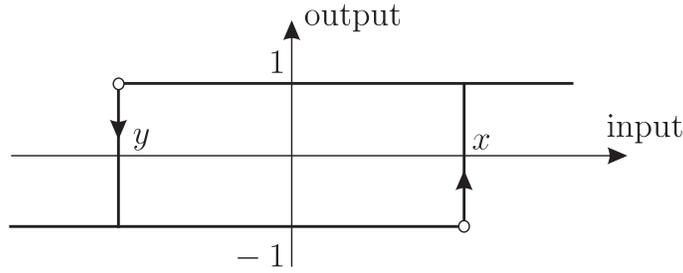}
\caption{Non-ideal relay.} \label{hyst}
\end{figure}

\item The population includes bacteria with different bi-stability
ranges $(-x,x)$, where the threshold value $x$ varies over an
interval $[\underline{x}, \overline{x}]\subset (0,1/2)$. We will
denote by $u(x,t)$ the density of the biomass of bacteria with
given switching thresholds $\pm x$ at a moment $t$.

\item There is a diffusion process acting on the density $u$.

\item At any particular time moment $t$, for any given $x$, all
the bacteria with the switching threshold values $\pm x$ are in
the same state (phenotype).
   That is, $u(x,t)$ is the total density of bacteria with the threshold $x$ at the moment $t$ and they are all in the same state.
This means that when a bacterium with a threshold $x'$
sporadically changes its threshold to a different value $x$, it
simultaneously copies the state from other bacteria which have the
threshold~$x$. In particular, this may require a bacterium to
change the state when its threshold changes.
\end{enumerate}

With these assumptions, we obtain the following model of the
evolution of bacteria and nutrients,
\begin{equation}\label{00}
\left\{
\begin{aligned}
& {u}_t=  D u_{xx} + \frac12  (1+{\mathcal R}^x(w))\, u
f_1+\frac12(1-{\mathcal R}^x(w))\, u f_{-1},\\
& \dot{f}_{1} =-\frac12
f_1\int_{\underline{x}}^{\overline{x}}(1+{\mathcal R}^x(w))\,
u\,dx,\\
& \dot{f}_{-1} =-\frac12
f_{-1}\int_{\underline{x}}^{\overline{x}}(1-{\mathcal R}^x(w))\,
u\,dx,
\end{aligned}
\right.
\end{equation}
where $u_t$ and $u_{xx}$ are the derivatives of the population
density $u$;\ $D>0$ is the diffusion coefficient; dot denotes the
derivative with respect to time; and all the non-ideal relays
${\mathcal R^x}$, $x\in[\underline{x},\overline{x}]$, have the
same input $w=f_1/(f_1+f_{-1})-1/2$. We additionally assume the
growth rate $\frac12 (1+ i {\mathcal R}^x(w))\, u f_i $ based on
the mass action law for bacteria in the phenotype $i=\pm 1$.  The
rate of the consumption of nutrient in the equation for
$f_i=f_i(t)$ is proportional to the total biomass of bacteria in
the phenotype $i$, hence to the integral; $\underline{x}$ and
$\overline{x}$ are the lower and upper bounds on available
threshold values, respectively.

We assume that a certain amount of nutrients is available at the
initial moment; the nutrients are not supplied after that moment.
 We assume the Neumann boundary conditions for $u$, that is no flux
of the population density $u$ through the lower and upper bounds
of available threshold values.

\section{Rigorous setting of a well-posed model}\label{secSetting}

\subsection{Rigorous setting}

Throughout the paper,  we assume  that
$x\in[\ux,\ox]\subset(0,1/2)$.

We begin with a rigorous definition of the hysteresis operator
$\cR^x$ ({\em non-ideal relay}) with fixed thresholds $\pm x$.
This  operator takes a continuous function $w=w(t)$ defined on an
interval $t\in[0,T)$ to the binary function $r=\cR^x(w)$ of time
defined on the same interval, which is given by
\begin{equation}\label{nonideal}
 {\mathcal R}^x (w)(t) = \left\{
\begin{array}{rll} -1 & {\rm if} & w(\tau)\le -x \ {\rm for\ some} \ \tau\in [0,t] \\
 && {\rm and} \ w(s)<x \ {\rm for\ all} \ s\in [\tau,t],
 \\1& {\rm if} &  w(\tau)\ge x \ {\rm for\ some} \ \tau\in [0,t] \\
 && {\rm and} \ w(s)>-x \ {\rm for\ all} \ s\in [\tau,t],\\
 r_0 & {\rm if} & -x < w(\tau) < x \ {\rm for\ all} \ \tau\in
 [0,t],
 \end{array}
\right.
\end{equation}
where $r_0$ is either $1$ or $-1$ (initial state of the non-ideal
relay $\cR^x$). Since $r_0$ may take different values for
different $x$, we write $r_0=r_0(x)$. The function $r_0=r_0(x)$ of
$x\in[\underline{x},\overline{x}]$ taking values $\pm 1$ is called
the {\em initial configuration} of the non-ideal relays. In what
follows, we  do not explicitly indicate the dependence of the
operator ${\mathcal R}^x$ on $r_0(x)$.

In this paper, we assume that $r_0(x)$ is {\it simple}, which
means the following. There is a partition $\underline{x}=\bar
x_0<\bar x_1<\cdots<\bar x_{N_0}=\overline{x}$ of the interval
$[\underline{x},\overline{x}]$ such that the function $r_0(x)$,
which satisfies $|r_0(x)|=1$ for all $\underline{x}\le x\le
\overline{x}$, is constant on each interval $(\bar x_{k-1},\bar
x_k]$ and has different signs on any two adjacent intervals:
\begin{equation}\label{simple}
\begin{array}{l}
r_0(x)=r_0(\bar x_k),\qquad x\in(\bar x_{k-1},\bar x_k],\ k=1,\ldots,N_0,\\
r_0(\bar x_{k-1})r_0(\bar x_k)=-1,\ k=2,\ldots,N_0,
\end{array}
\end{equation}
where the second relation holds if $N_0\ge2$.

We define the {\em distributed relay operator} $\cR(w)$ taking
functions $w=w(t)$ to functions $r=r(x,t)$ by
\begin{equation}\label{eqrxt}
r(x,t)=\cR(w)(x,t):={\mathcal R}^x(w)(t).
\end{equation}
The function $r(\cdot,t)$ will be referred to as the configuration
(state) of the distributed relay operator at the moment $t$.

We set
\begin{equation}\label{eqNP}
\uU(t)= \int_{\underline{ x}}^{\overline{ x}} u( x,t)\,d x,\quad
\cP(u,w)(t)=\int_{\underline{ x}}^{\overline{ x}} u( x,t)
\cR^x(w)(t)\,d x.
\end{equation}
Here the first integral is the total mass of bacteria; $\cP$ is
the so-called {\em Preisach operator} with the time dependent
density function $u$. Further, we replace the unknown functions
$f_1$ and $f_{-1}$ in system \eqref{00} by the functions $v=f_1+f_{-1}$ (total
mass of the two nutrients) and $w=f_1/(f_1+f_{-1})-1/2$ (deviation
of the relative concentration of the first nutrient from the
value~1/2). Substituting
$$
f_1=\left(\dfrac{1}{2}+w\right)v,\qquad f_{-1}=\left(\dfrac{1}{2}-w\right)v
$$
into the first equation of \eqref{00},
we obtain the equation
$$
u_t = D u_{xx} +\left(\frac12 + w  \cR(w)\right)u v.
$$
Furthermore, summing the second and the third equations of system \eqref{00}
and using the relationships $v=f_1+f_{-1}$, $2 w v  =f_1-f_{-1}$ and the notation \eqref{eqNP}, we get
$$
\dot v=-\left(\frac{\uU}{2}  +w \cP(u,w)\right)v.
$$
Finally, from the second and the third equations of \eqref{00}, it follows that
$$
\dot w =\frac{\dot f_1 f_{-1} - f_1 \dot f _{-1}}{(f_1+f_{-1})^2}=-\frac{f_1 f_{-1}}{v^2}\int_{\underline{ x}}^{\overline{ x}} u \cR^x(w)\,d x=- \left(\dfrac{1}{2}+w\right)\left(\dfrac{1}{2}-w\right)
\cP(u,w).
$$

Thus,  the resulting system, which is equivalent to equations \eqref{00}, takes the form
\begin{equation}\label{eqBacteriaGeneral}\left\{
\begin{aligned}
& u_t = D u_{xx} +\left(\frac12 + w  \cR(w)\right)u v,\\
&\dot v=-\left(\frac{\uU}{2}  +w
\cP(u,w)\right)v,\\
&\dot w=- \left(\dfrac{1}{2}+w\right)\left(\dfrac{1}{2}-w\right)
\cP(u,w),
\end{aligned}\right.
\end{equation}
where we assume the Neumann boundary conditions
\begin{equation}
\label{eqBC}  u_x|_{ x=\ux} = u_x|_{ x=\ox}=0
\end{equation}
and the initial conditions
\begin{equation}
\label{eqIC}
\begin{gathered}
u( x,0)=u_0( x),\quad v(0)=v_0,\quad w(0)=w_0,\quad
r(x,0)=r_0(x),\\
u_0( x)\ge0,\quad \int_\ux^\ox u_0(x)\,dx =1,\quad v_0\ge 0, \quad
|w_0|\le \ox,\quad r(x)\ \text{is simple}.
\end{gathered}
\end{equation}

\subsection{Well-posedness}\label{secWellPosedness}
The problem  \eqref{eqBacteriaGeneral}--\eqref{eqIC}, which
contains a discontinuous distributed relay operator, was shown  in
\cite{GurRachDiffThresh2013} to be well posed. We briefly
summarize this result before proceeding with the analysis of long
time behavior.

Set  $Q_T=(\ux,\ox)\times(0,T)$ for $T>0$. We will use the
standard Lebesgue spaces $L_2(Q_T)$ and $L_ 2=L_ 2(\ux,\ox)$; the
Sobolev spaces $W_ 2^k=W_ 2^k(\ux,\ox)$, $k\in\bbN$; the
anisotropic Sobolev space $W_2^{2,1}(Q_T)$  with the norm
$$
\|u\|_{W_2^{2,1}(Q_T)}=\left(\int\limits_0^T
\|u(\cdot,t)\|_{W_2^2}^2\, dt+ \int\limits_0^T
\|u_t(\cdot,t)\|_{L_2}^2\, dt\right)^{1/2};
$$
and the space
$$
\cW(Q_T)=W_2^{2,1}(Q_T)\times C^1[0,T]\times C^1[0,T].
$$

\begin{definition}
Assume that $(u_0,v_0,w_0) \in W_2^1\times\bbR^2$. We say that
$(u,v,w)$ is a ({\em global})  {\em solution} to
problem~\eqref{eqBacteriaGeneral}--\eqref{eqIC}  if, for any
$T>0$,
 $(u,v,w)\in
\cW(Q_T)$, $r(\cdot,t)$ is a continuous $L_2$-valued function for
$t\ge0$, and relations~\eqref{eqBacteriaGeneral}--\eqref{eqIC}
hold in the corresponding function spaces.
\end{definition}


The following result was proved in \cite{GurRachDiffThresh2013}.

\begin{theorem}\label{thWellPosed}
If  $(u_0,v_0,w_0) \in W_2^1\times\bbR^2$, then
\begin{enumerate}
\item Problem~$\eqref{eqBacteriaGeneral}$--$\eqref{eqIC}$ has a
unique solution $(u,v,w)$$;$

\item The state $r(\cdot,t)=\cR(w)(\cdot,t)$ of the distributed
relay operator is {simple} for all $t\ge 0$$;$

\item We have
\begin{gather} \dot \uU(t)\ge0,\quad \uU(t)\to
1+v_0,\quad
 0\le u(\cdot,t)\to\dfrac{1+v_0}{\ox-\ux}\ \ \text{in}\ \
 C[\ux,\ox], \label{limit''}\\
\dot v(t)\le 0,\qquad v(t)\le v_0e^{-\mu t} \quad \text{with} \quad \mu=1/2-\overline{x}>0, \label{limit'} \\
 |w(t)|\le\ox<1/2, \label{exist}
\end{gather}
where   the convergence takes place as $t\to\infty$.
\end{enumerate}
\end{theorem}

The   behavior given by~\eqref{limit''}--\eqref{exist} is to be
expected. Indeed, as we assume no supply of nutrients after the
initial moment, the total amount of nutrients $v(t)$ converges to
zero. When the density of nutrients vanishes as a result of
consumption by bacteria, the equation for the density $u$
approaches the homogeneous heat equation with zero flux boundary
conditions, which explains why the density of bacteria $u(x,t)$
converges to a uniform distribution over the interval
$[\underline{x},\overline{x}]$ as a result of the diffusion.

Below we are interested in the limiting behavior of the
configuration function $r(x,t)$ given by~\eqref{eqrxt}, or, in
other words, in the limiting distribution of phenotypes over the
population of thresholds. In~\cite{GurRachPattern2015}, we have
proved that $r(x,t)$ converges to a step like profile as
$t\to\infty$, see Fig.~\ref{figSteps}. In the next section, we give a quantative
description of this phenomenon under the assumptions that the
initial amount of nutrient $v_0$ and the diffusion coefficient $D$
tend to zero, while the initial density $u_0(x)$ tends to the
delta function.






\section{Fronts asymptotics}\label{secFrontAsymp}

\subsection{Terminology}

  A point $\bar x_j=\bar x_j(t)$
which separates an interval on which the relays are in state $1$
from  an interval on which the relays are in state $-1$ is called
a {\em front}; cf.~\eqref{simple}. The total number of fronts can
vary, but, by Theorem~\ref{thWellPosed}, remains finite at all
times (equivalently, the state of the distributed relay operator
remains simple at all times). A front can either stay (a {\em
steady} front) or move right. That is, any $\bar x_j(t)$ is a
non-decreasing function on the time interval of its existence. A
front {\em disappears} if it hits another front or is hit by
another front. Assume a front $\bar x_j(t)$ exists at a moment
$t_0$. It is called {\em immortal up to time $T$}, $T>t_0$, if it
does not disappear on the time interval $t_0\le t\le T$.

Our goal is to give asymptotic formulas for the time moments at
which steady immortal (up to time $T$) fronts appear as well as
for their positions. This will be done as $D\to 0$, where $D$ is the diffusion coefficient in~\eqref{eqBacteriaGeneral} responsible for   random fluctuations of hysteresis thresholds of individual bacteria.

\subsection{Assumptions}\label{subsecAssump}

From now on, we fix a time interval $[0,T]$, $T>0$. Along with
conditions~\eqref{eqIC}, we assume that the initial density
satisfies the relation
\begin{equation}\label{data2}
u_0(x)\le\varepsilon,\quad  x\in
(\underline{x},\overline{x}-\varepsilon),
\end{equation}
where $\varepsilon>0$ will be chosen small enough and depending on
$D$.

We introduce the function
\begin{equation}\label{eqFx}
 F(x):={1}/{4}-(\ox-x)^2
\end{equation}
  and set
\begin{equation}\label{eqs12}
s_{1/2}:=\int_0^{\ox+\ux}\dfrac{dx}{F(x)},
\end{equation}
\begin{equation}\label{eqoverF}
\oF:=2\int_0^\ox\dfrac{dx}{F(x)}.
\end{equation}
Set
$$
U(x,t):=\int_x^{\overline x} u(y,t)\,dy,\qquad
E(y):=\frac{2}{\sqrt{\pi}}\int_0^y e^{-z^2}\,dz.
$$

Using~\cite[Lemmas~4.1 and~5.4]{GurRachPattern2015} and the fact
that $0\le\uU(t)-1\le v_0$ (see~\eqref{eqIC} and~\eqref{limit''}),
we can choose positive functions $\varepsilon(D)$ and $v_0(D)$
such that $\varepsilon(D),v_0(D)\to0$ as $D\to 0$ and
\begin{equation}\label{eqUcloseE}
\sup\limits_{t\in[s_{1/2}/2,T]}\left\|U(\cdot,t)-
E\left(\dfrac{\ox-\cdot}{2(D
t)^{1/2}}\right)\right\|_{C[\ux,\ox]}=o(1)\quad \text{as } D\to 0,
\end{equation}
provided that $u_0(x)$ satisfies~\eqref{data2} with
$\varepsilon=\varepsilon(D)$ and that $v_0=v_0(D)$.

The main assumptions will be as follows.
\begin{enumerate}
\item The initial density $u_0(x)$ satisfies~\eqref{data2} with
$\varepsilon=\varepsilon(D)$, where $\varepsilon(D)$ is the above
function.
%

\item The initial amount of food is small: $v_0=v_0(D)$, where
$v_0(D)$ is the above function.

\item The initial value of the input is close to $\ox$:
$w_0=w_0(D)$, where $w_0(D)$ is an arbitrary function such that
$w_0(D)\le\ox$ and $w_0(D)\to\ox$ as $D\to 0$.

\item The initial configuration is $r_0(x)\equiv 1$.
\end{enumerate}


Our goal is to determine  the consecutive time moments $t_n$  at
which the moving fronts become steady and their positions $x_n$
($n=1,2,\dots$)   at these moments.

We begin with the following recursive algorithm for determining
(finitely or infinitely many) positive numbers
\begin{equation}\label{eqSequencesnyn}
0=s_0<s_1<s_2<\dots<s_n<\dots,\quad  y_1,y_2,\dots,y_n,\dots
\end{equation}

\subsection{Algorithm}\label{subsecAlgorithm}

{\bf Basis.} Set $s_1:=\oF$, and let $y_1$ be the $($unique$)$
root of the equation $2E(y)=1$, $y\in(0,\infty)$.

Set $G_2(t):=-2
E\left(y_1\left(\dfrac{s_1}{t}\right)^{1/2}\right)$. Note that
$G_2(s_1)+1=0$ and $\dot G_2(t)>0$ for all $t>0$.

{\bf Inductive conjecture.} Fix $n\ge 2$. Assume that we have defined
the sequences $s_1,\dots,s_{n-1}$ and $y_1,\dots,y_{n-1}$ such that the function
\begin{equation}\label{eqGn}
G_n(t):= 2\sum\limits_{j=1}^{n-1}(-1)^{n+j}
E\left(y_j\left(\dfrac{s_j}{t}\right)^{1/2}\right)
\end{equation}
satisfies
\begin{equation}\label{eqGnsn1}
G_n(s_{n-1})+1=0,\qquad   G_n(t)+1>0,\quad t>s_{n-1}.
\end{equation}
As we have seen, relations~\eqref{eqGnsn1} hold for $n=2$.

Below we will also use the equivalent recursive representation of
$G_n(t)$:
\begin{equation}\label{eqGnRecursive}
G_1(t)\equiv 0,\quad G_n(t)=-G_{n-1}(t)-2
E\left(y_{n-1}\left(\dfrac{s_{n-1}}{t}\right)^{1/2}\right),\quad
n\ge 2.
\end{equation}

{\bf Inductive step.} Now we will determine $s_n$, $y_n$, and the
function $G_{n+1}(t)$ satisfying relations~\eqref{eqGnRecursive}
with $n$ replaced by $n+1$.
\begin{enumerate}

\item\label{caseAExistsTildeS} Consider the equation
\begin{equation}\label{eqsn}
\int_{s_{n-1}}^{s}(G_n(t)+1)\,dt=\oF,\quad s\in[s_{n-1},\infty).
\end{equation}
Due to~\eqref{eqGnsn1} and the fact that $G_n(t)+1\to 1$ as
$t\to\infty$, the left-hand side is a strictly increasing
function of $s$ that tends to $\infty$ as $s\to\infty$. Therefore,
Eq.~\eqref{eqsn} has a unique root, which we denote by  $\tilde
s_n$.

\item Consider the function
\begin{equation}\label{eqHtz}
H_n(t,z):=-G_n(t)-2E\left(\dfrac{z}{t^{1/2}}\right)+1,\quad t\ge
\tilde s_n,\ z\ge 0.
\end{equation}
Relations~\eqref{eqGnsn1} imply that
\begin{equation}\label{eq-Gn12}
-G_n(t)+1< 2,\quad t\ge \tilde s_n.
\end{equation}
 On the other hand, due
to~\eqref{eqGnRecursive},
$$
-G_n(t)+1=G_{n-1}(t)+2
E\left(y_{n-1}\left(\dfrac{s_{n-1}}{t}\right)^{1/2}\right)+1,\quad
t\ge \tilde s_n,
$$
and, hence,
\begin{equation}\label{eq-Gn10}
-G_n(t)+1>0,\quad t\ge \tilde s_n,
\end{equation}
 due to the inductive conjecture (see~\eqref{eqGnsn1}
with $n$ replaced by $n-1$).

Thus, we see that the equality $H_n(t,z)=0$ uniquely defines a
smooth {\it positive} function $z=Z_n(t)$, $t\ge \tilde s_n$, such
that
\begin{equation}\label{eqZn}
H_n(t,z)>0,\quad z<Z_n(t);\qquad H(t,Z_n(t))=0;\qquad
H_n(t,z)<0,\quad z>Z_n(t).
\end{equation}


\item Now we have several possibilities depending on the geometry
of the curve $Z_n(t)$ defined by~\eqref{eqZn}:
\begin{enumerate}
\item\label{caseA}   $\dot Z_n(t)>0$ for $t>\tilde s_n$. In this
case, we set $s_n:=\tilde s_n$.

\item\label{caseB} There is $s_n>\tilde s_n$ such that $\dot
Z_n(t)<0$ for $t\in(\tilde s_n,s_n)$ and $\dot Z_n(t)>0$ for
$t>s_n$.

\item Otherwise, we terminate the process at the finite  sequences
$\{s_1,\dots,s_{n-1}\}$ and $\{y_1,\dots,y_{n-1}\}$.
\end{enumerate}

In cases~\ref{caseA} and~\ref{caseB}, we set
\begin{equation}\label{eqFormulayn}
y_n:=\dfrac{Z_n(s_n)}{s_n^{1/2}}.
\end{equation}

In these two cases, it remains to check that the function
\begin{equation}\label{eqGnPlus1}
G_{n+1}(t):=-G_{n}(t)-2
E\left(y_{n}\left(\dfrac{s_{n}}{t}\right)^{1/2}\right)
\end{equation}
satisfies relations~\eqref{eqGnsn1} with $n$ replaced by $n+1$.
This will guarantee that we can repeat the inductive step of the algorithm from section~\ref{subsecAlgorithm} with $n$
replaced by $n+1$.

Combining~\eqref{eqGnPlus1} with~\eqref{eqHtz}, \eqref{eqZn}, and
\eqref{eqFormulayn},    we have
$$
G_{n+1}(s_n)+1=-G_{n}(s_n)-2
E\left(y_{n}\right)+1=H(s_n,Z(s_n))=0.
$$
On the other hand, in cases~\ref{caseA} and~\ref{caseB}, we have
$Z_n(s_n)<Z_n(t)$ for $t>s_n$. Therefore, relations~\eqref{eqHtz},
\eqref{eqFormulayn}, and~\eqref{eqGnPlus1} imply that
$$
G_{n+1}(t)+1=H_n(t,Z_n(s_n))>0,\quad t>s_n.
$$
\end{enumerate}

\subsection{Main result}

\subsubsection{Formulation of the main result}

Assume that, for some $N\ge 1$, we have the sequences
$\{s_1,\dots,s_{N}\}$ and $\{y_1,\dots,y_{N}\}$ constructed
according to the above algorithm, and assume that $s_N< T$, where
$T$ was fixed in Sec.~\ref{subsecAssump}.

\begin{theorem}\label{tAsympImmortalFronts}  For all $n=1,\dots,N$, the
  $n$-th moving front becomes
steady and immortal $($up to time $s_N)$ at a moment
\begin{equation}\label{eqtnsno1}
t_n=s_n+o(1)
\end{equation}
and its position at this moment is
\begin{equation}\label{eqxnqn}
x_n= \ox-q_n ,\quad\text{where}\quad
q_n=2\big(y_n+o(1)\big)\big(Ds_n\big)^{1/2};
\end{equation}
here $o(1)$ stands for functions of   $D$ that tend to $0$ as
$D\to 0$.
Furthermore, $x_N<x_{N-1}<\dots<x_1<\ox$.
\end{theorem}

\subsubsection{Preliminary discussion}\label{subsecPreliminaryDiscussion}

In section~\ref{secProofMainResult}, we will prove
Theorem~\ref{tAsympImmortalFronts}. Assume we have constructed
$n-1$ fronts that became steady at the moments $t_1,\dots,t_{n-1}$
at the positions $x_1,\dots,x_{n-1}$ given by~\eqref{eqtnsno1}
and~\eqref{eqxnqn}, respectively. We will consider the formation
of the $n$-th front on the time interval $[t_{n-1},t_n]$. We
consider the case of an even $n\ge 2$. The case of odd $n$ is analogous;
see also Remark~\ref{remn1} below concerning $n=1$.

Below we will consecutively consider time moments
$$
\begin{aligned}
 \big[t_{n-1}& =s_{n-1}+o(1)\big]<r_{n-1}<t_{n-1/2}\\
& <\big[\tilde\tau_n=\tilde s_n+o(1)\big] \le \big[\tau_n=s_n+o(1)\big]\\
& \le\big[\theta_n=s_n+o(1)\big]  \le \big[t_n=s_n+o(1)\big]<r_n
\end{aligned}
$$
which are characterized by the following properties:
\begin{enumerate}
\item The function $w(t)$ is {\em increasing} and the relays {\em
do not switch} during the time interval $t\in(t_{n-1},t_{n-1/2})$; furthermore,
\begin{equation}\label{eqtimes1}
w(t_{n-1})=-\ox+q_{n-1},\quad w(t_{n-1/2})=\ux.
\end{equation}

\item The function $w(t)$ is {\em increasing} and some relays {\em
switch} during the time interval $(t_{n-1/2},\tilde\tau_n)$; further, for some
$Z_*(D)=o(1)$ (which will be chosen below to satisfy relations~\eqref{eqHnzeta}), we have
\begin{equation}\label{eqtimes2}
w(\tilde\tau_n)=\ox-2Z_*(D).
\end{equation}
Note that, although the distance between $w(\tilde\tau_n)$ and $\ox$ is already $2Z_*(D)=o(1)$, it is still much larger than the desired distance $2\big(y_n+o(1)\big)\big(Ds_n\big)^{1/2}$ due to~\eqref{eqZ*OD12}. Also the moment~$\tilde\tau_n$, which is of order $\tilde s_n$, may be much smaller than the desired moment, which is of order $s_n$.

\item The function $w(t)$ {\em may  oscillate} and the relays {\em
may, though not necessarily, switch} during the time interval $(\tilde\tau_n,\tau_n)$.
 The moment $\tau_n=s_n+o(1)$ is already of the correct order, but the value $w(\tau_n)$ may still have a wrong asymptotics. However, it satisfies
\begin{equation}\label{eqtimes2'}
w(\tilde\tau_n) \le
w(\tau_n)\le\ox-2\big(y_n+o(1)\big)\big(Ds_n\big)^{1/2},
\end{equation}
i.e., it is bounded from above by the correct asymptotics.

  \item The function $w(t)$ {\em may  oscillate} and the relays {\em
may, though not necessarily, switch} during the time interval $(\tau_n,\theta_n)$. The moment $\theta_n=s_n+o(1)$   is of the same (correct) order as~$\tau_n$, but now, additionally, the value $w(\theta_n)$ has the correct asymptotics
\begin{equation}\label{eqtimes2''}
w(\theta_n)=\ox-2\big(y_n+o(1)\big)\big(Ds_n\big)^{1/2}.
\end{equation}

\item The function $w(t)$ {\em may  oscillate} and the relays {\em
may, though not necessarily, switch} during the time interval $(\theta_n,t_n)$. The moment $t_n=s_n+o(1)$ is of the same (correct) order as $\tau_n$ and $\theta_n$, and the function $w(t)$ does not oscillate any more after the moment $t_n$. More precisely, $w(t)$ is {\em decreasing} and the relays {\em
do not switch} during the time interval $(t_n,r_n]$. Moreover, the $n$-th immortal front is formed at the moment $t_n$ in the sense that
\begin{equation}\label{eqtimes3}
 w(t_n)=\ox-q_n=x_n,\quad \dot w(t_n)=0.
\end{equation}

\end{enumerate}


\section{Proof of the main result}\label{secProofMainResult}

In this section, we will prove Theorem~\ref{tAsympImmortalFronts}. Recall that, without loss of generality, we assume that $n$ is even. In sections~\ref{subsubsec1}--\ref{subsubsec5} below, we will consider in detail the respective time intervals from items 1--5 of section~\ref{subsecPreliminaryDiscussion}.

\subsection{Dynamics for $w\in (-\ox+q_{n-1},\ux)$ as $t\in(t_{n-1},t_{n-1/2})$}\label{subsubsec1}

  We will
see below that $w(t)$ is increasing for $t> t_{n-1}$ and achieves
the value $w(t)=\ux$. Let $t_{n-1/2}$ be the  moment when this
happens: $w(t_{n-1/2})=\ux$. Note that if $w(t)\in
(-\ox+q_{n-1},\ux)$ and $\dot w(t)>0$, then the relays do not
switch.

\begin{figure}[ht]
  \centering
\includegraphics[width=155mm]{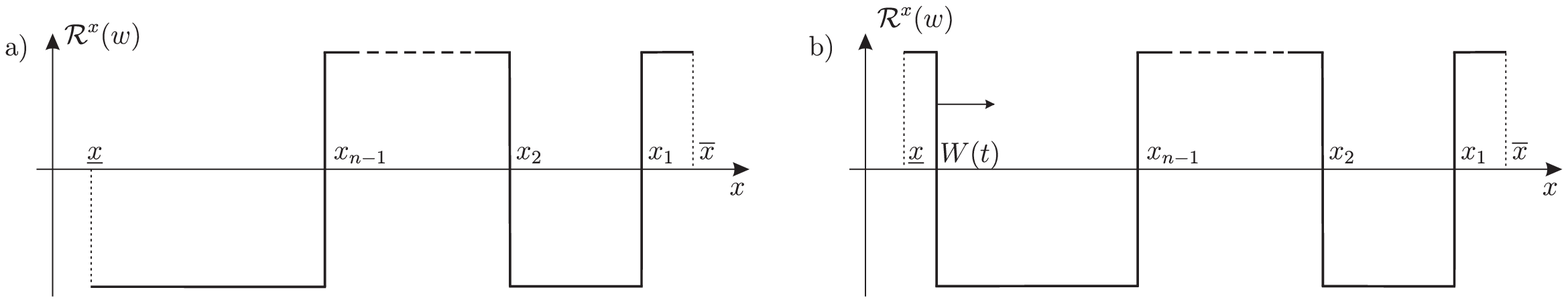}
\caption{Spatial configuration of the relays $\mathcal R^x(w)$ at
a moment $t>t_{n-1}$. a) The input satisfies $w(t)\in
(-\ox+q_{n-1},\ux)$. Hence, the relays do not switch and all the
fronts are steady. b) The input satisfies $w(t)\in (\ux,\ox-q_n)$.
Hence, the relays may switch and the $n$-th front $W(t)$ can move.} \label{figSteps}
\end{figure}

 Assume we have proved at the previous
step (for $n$ replaced by $n-1$) that
\begin{equation}\label{eqdotqrn1}
\dot w(t)>0 \quad\text{for } t\in(t_{n-1},r_{n-1}]
\end{equation}
 for some $r_{n-1}>s_{n-1}$
not depending on $D$ (we shall do this for $k=n$ in the end;
see~\eqref{eqdotqn0tn}). As long as $\dot w(t)$ remains positive
and $w(t)<\underline{x}$ holds, the relays do not switch, see
Fig~\ref{figSteps}.a. Hence, using~\eqref{eqUcloseE}
and~\eqref{eqxnqn}, we can write the third equation
in~\eqref{eqBacteriaGeneral} as follows (we omit ``$u(x,t)\,dx$''
for brevity):
\begin{equation}\label{eqdqdtn-1/2}
\begin{aligned}
\dot w(t)&=-\left(\dfrac{1}{4}-w^2\right)\left(-\int_\ux^{x_{n-1}}
+ \int_{x_{n-1}}^{x_{n-2}}   - \int_{x_{n-2}}^{x_{n-3}}
+\dots+\int_{x_{1}}^\ox  \right)\\
&=
\left(\dfrac{1}{4}-w^2\right)\big(U(\ux,t)-2U(x_{n-1},t)+2U(x_{n-2},t)-\dots-2U(x_1,t)\big)\\
& =
\left(\dfrac{1}{4}-w^2\right)\left(1+2\sum\limits_{j=1}^{n-1}(-1)^j
E\left((y_j+o(1))\left(\dfrac{s_j}{t}\right)^{1/2}\right)+o(1)\right),
\quad t>t_{n-1}
\end{aligned}
\end{equation}
(here and below, we consider $t\le T$ only.)

 We set
$$p(t):=\ox+w(t).$$ Taking into account~\eqref{eqdqdtn-1/2},
definition~\eqref{eqFx} of $F(x)$, definition~\eqref{eqGn} of
$G_n(t)$, and the fact that $n$ is even, we see that $p(t)$
satisfies
\begin{equation}\label{eqdqdtn-1/2Phi}
\dot p(t)=F(p)\big( G_n(t)+1+\mu(t,D)\big),\quad t>t_{n-1},
\end{equation}
\begin{equation}\label{eqqn-1/2_0}
p(t_{n-1})=2(y_{n-1}+o(1))(D s_{n-1})^{1/2},
\end{equation}
where
\begin{equation}\label{eqepsn}
\sup\limits_{t\in[s_{1/2}/2,T]}|\mu(t,D)|\le
\varepsilon_n=\varepsilon_n(D)\quad\text{for a nonnegative}\
\varepsilon_n(D)=o(1).
\end{equation}

Due to~\eqref{eqGnsn1} and~\eqref{eqepsn}, the right-hand side
in~\eqref{eqdqdtn-1/2Phi} is positive for $t\in[r_{n-1},T]$,
provided $D$ is small enough. Together with~\eqref{eqdotqrn1},
this means that the solution to~\eqref{eqdqdtn-1/2Phi},
\eqref{eqqn-1/2_0} increases until a moment $t_{n-1/2}$ (if it
exists) at which it achieves the value $\ox+\ux$. This proves
that~\eqref{eqdqdtn-1/2Phi} is  equivalent to the third equation
in~\eqref{eqBacteriaGeneral} not only for $t\in(t_{n-1},r_{n-1})$
but actually for $t\in(t_{n-1},t_{n-1/2})$.

The following lemma determines the time moment $t_{n-1/2}$.

\begin{lemma}\label{tn-1/2}
The equation
\begin{equation}\label{eqsn-1/2'}
\int_{s_{n-1}}^{s}(G_n(t)+1)\,dt=\int_0^{\ox+\ux}\dfrac{dx}{F(x)}
\end{equation}
has a unique root on the interval $s\in(s_{n-1},\infty)$. Denote
it by $s_{n-1/2}$. Then $t_{n-1/2}=s_{n-1/2}+o(1)$ and $\dot
p(t_{n-1/2})>0$.
\end{lemma}
\proof Due to~\eqref{eqGnsn1} and the fact that $G_n(t)+1\to 1$ as
$t\to\infty$, Eq.~\eqref{eqsn-1/2'} has a unique root $s_{n-1/2}$.

To prove the asymptotics for $t_{n-1/2}$, we denote by
$p_{\pm}(t)$ the solutions of
$$
\dot p_{\pm}(t)=F(p_\pm)\big( G_n(t)+1\pm\varepsilon_n\big),\quad
t>t_{n-1},
$$
with the same initial data at $t=t_{n-1}$ as
in~\eqref{eqqn-1/2_0}. Obviously,
\begin{equation}\label{eqsn-1/2''}
p_-(t)\le p(t)\le p_+(t).
\end{equation}

Let $t_{n-1/2,\pm}$ be the first moment at which
$p_\pm(t)=\ox+\ux$. Then, integrating the   differential equation
for $p_\pm(t)$, we obtain
$$
\int_{s_{n-1}+o(1)}^{t_{n-1/2,\pm}}(G_n(t)+1\pm\varepsilon_n)\,dt=\int_{o(1)}^{\ox+\ux}\dfrac{dx}{F(x)}.
$$
Using the implicit function theorem in a neighborhood of the point
$(t_{n-1/2,\pm},D)=(s_{n-1/2},0)$, we obtain
$t_{n-1/2,\pm}=s_{n-1/2}+o(1)$. Therefore, taking into
account~\eqref{eqsn-1/2''}, we have $t_{n-1/2}=s_{n-1/2}+o(1)$.

The inequality  $\dot p(t_{n-1/2})>0$ has been proved before the
lemma.
\endproof

\subsection{Dynamics for $w\in (\ux,\ox-q_n)$ as $t\in(t_{n-1/2},\tilde\tau_n)$}\label{subsubsec2}

By Lemma~\ref{tn-1/2}, $\dot w(t_{n-1/2})>0$. Therefore, as long
as $\dot w(t)$ remains positive, the input $w(t)$ increases and
hence switches the relays. However, if $\dot w(t)$
becomes negative, the relays will stop switching. To describe the
dynamics of $w(t)$, we introduce the function
$$
W(t):=\max\limits_{s\in[t_{n-1/2},t]}w(s),
$$
where $t_{n-1/2}$ is defined in Lemma~\ref{tn-1/2}, see
Fig~\ref{figSteps}.b. Then, similarly to~\eqref{eqdqdtn-1/2}, we
obtain from the third equation in~\eqref{eqBacteriaGeneral} and
from relation~\eqref{eqUcloseE}
\begin{equation*}
\begin{aligned}
\dot
w(t)&=-\left(\dfrac{1}{4}-w^2\right)\left(\int_\ux^{W(t)}-\int_{W(t)}^{x_{n-1}}
+ \int_{x_{n-1}}^{x_{n-2}}   - \int_{x_{n-2}}^{x_{n-3}}
+\dots+\int_{x_{1}}^\ox  \right)\\
&=
-\left(\dfrac{1}{4}-w^2\right)\big(1-2U(W,t)+2U(x_{n-1},t)-2U(x_{n-2},t)+\dots+2U(x_1,t)\big)\\
& =
-\left(\dfrac{1}{4}-w^2\right)\left(1-2E\left(\dfrac{\ox-W}{2(D
t)^{1/2}}\right)-2\sum\limits_{j=1}^{n-1}(-1)^j
E\left((y_j+o(1))\left(\dfrac{s_j}{t}\right)^{1/2}\right)+o(1)\right),
\quad t>t_{n-1/2}.
\end{aligned}
\end{equation*}
 Setting
$$
q(t):=\ox-w(t),\quad
Q(t):=\min\limits_{s\in[t_{n-1/2},t]}q(s)=\ox-W(t),
$$
we have
\begin{equation}\label{eqdqdtnPhiMin}
\dot q(t)=F(q)\left(-G_n(t)- 2E\left(\dfrac{Q}{2(D
t)^{1/2}}\right)+1+\mu(t,D)\right),\quad t>t_{n-1/2},
\end{equation}
 where
\begin{equation}\label{eqMuoD12}
 \sup\limits_{t\in[s_{1/2}/2,T]}|\mu(t,D)|\to
0\quad\text{as}\ D\to 0
\end{equation}
with $s_{1/2}$ defined in~\eqref{eqs12}.
 We consider Eq.~\eqref{eqdqdtnPhiMin} with the
initial condition
\begin{equation}\label{eqqn_0}
q(t_{n-1/2})=\ox-\ux.
\end{equation}
We begin with the following observation, which shows that the
$(n-1)$-th front is immortal.
\begin{lemma}\label{lwnwn-1}
The solution $q(t)$ to problem~\eqref{eqdqdtnPhiMin},
\eqref{eqqn_0} cannot achieve the value $q_{n-1}$.
\end{lemma}
\proof Assume, by contradiction, that $q(t)$ achieves the value $q_{n-1}$
($=2(y_{n-1}+o(1))(D s_{n-1})^{1/2})$) for the first time at some
moment $t=t^*$ ($>t_{n-1/2}>s_{n-2}$). At this moment
$Q(t^*)=q(t^*)=q_{n-1}$. Since $t^*>s_{n-2}$, it follows
from~\eqref{eqGnsn1} (with $n$ replaced by $n-1$) that $
G_{n-1}(t^*)+1>0$. Therefore, using~\eqref{eqdqdtnPhiMin} and
\eqref{eqGnRecursive}, we have
$$
\dot q(t^*)=F(0)\big(G_{n-1}(t^*)+1\big)+o(1)>0
$$
for all sufficiently small $D$, which is  impossible.
\endproof

%


To proceed, we will need a specific (positive) function $Z_*(D)$
such that
\begin{equation}\label{eqZ*o1}
Z_*(D)\to 0\quad\text{as }D\to 0,
\end{equation}
\begin{equation}\label{eqZ*OD12}
\dfrac{Z_*(D)}{D^{1/2}}\to \infty\quad\text{as }D\to 0.
\end{equation}
 To define it, we introduce the new function
$z(t):=\dfrac{q(t)}{2D^{1/2}}$ and set
\begin{equation}\label{eqZminz}
Z(t):=\min\limits_{s\in[t_{n-1/2},t]}z(s).
\end{equation}
Then we consider  Eq.~\eqref{eqdqdtnPhiMin} for $t>\tilde s_n$,
where $\tilde s_n$ is the time moment defined in the inductive
step~\ref{caseAExistsTildeS} of the algorithm from section~\ref{subsecAlgorithm}. For such $t$, it is equivalent to
\begin{equation}\label{eqdzdt2Phi}
 \dot z(t) = \dfrac{F(0)(H_n(t,Z)+o(1))}{2D^{1/2}},
\end{equation}
where $H(t,\cdot)$ is given by~\eqref{eqHtz}. We choose the
desired  $Z_*(D)$ to satisfy, for all $t\in[\tilde s_n,T]$,
\begin{equation}\label{eqHnzeta}
\begin{aligned}
& H_n(t,\zeta)+o(1)>0 & & \text{for } \zeta\le Z_n(t)-Z_*(D),\\
& H_n(t,\zeta)+o(1)<0 & & \text{for } \zeta\ge Z_n(t)+Z_*(D),\\
& F(0)|H_n(t,\zeta)+o(1)|\le 2 h_n Z_*(D) & & \text{for } \zeta\in[Z_n(t)-Z_*(D),Z_n(t)+Z_*(D)],\\
\end{aligned}
\end{equation}
where $o(1)$ is the same function as in~\eqref{eqdzdt2Phi} and
$h_n>0$ does not depend on $t$ and $D$. This can be done because
$\partial H_n(t,Z_n(t))/\partial\zeta\ne 0$ for all $t\in[\tilde
s_n,T]$.

 The following lemma defines a time moment
$\tilde\tau_n$ satisfying~\eqref{eqtimes2}.

\begin{lemma}\label{ltau2s2}
There is a time moment $ \tilde\tau_n=\tilde s_n+o(1)$ such that
$\dot q(t)<0$ for all $t\in[t_{n-1/2},\tilde\tau_n]$ and
$q(\tilde\tau_n)=2 Z_*(D)$.
\end{lemma}
\proof As long as $\dot q(t)<0$ and $Z_*(D)\le q(t)\le \ox-\ux$,
we have $Q(t)=q(t)$. Therefore, we can rewrite
Eq.~\eqref{eqdqdtnPhiMin} as follows:
\begin{equation}\label{eqdqdtnPhitau2s2_1}
\dot q(t)=-F(q)\big(G_n(t)+1+\mu_1(t,q,D)\big),\quad t>t_{n-1/2},
\end{equation}
where, due to~\eqref{eqMuoD12} and~\eqref{eqZ*OD12},
\begin{equation}\label{eqdeltaDo1}
\sup\limits_{q\in[0,\ox-\ux]}\sup\limits_{t\in[s_{1/2}/2,T]}|\mu_1(t,q,D)|\le
\delta=\delta(D) \quad\text{for a nonnegative}\ \delta(D)=o(1).
\end{equation}

Due to~\eqref{eqGnsn1} and~\eqref{eqdeltaDo1}, the right-hand side
in~\eqref{eqdqdtnPhitau2s2_1} is negative for $t\in[t_{n-1/2},T]$,
provided $D$ is small enough. Hence, $\dot q(t)<0$ for
$t>t_{n-1/2}$.

Denote by $q_{\pm}(t)$ the solution to the equation
\begin{equation}\label{eqqpmeps}
\dot q_{\pm}(t)=-F(q_\pm)\big(G_n(t)+1\pm\delta\big),\quad t>
t_{n-1/2},
\end{equation}
with the same initial condition as in~\eqref{eqqn_0},
$q(t_{n-1/2})=\ox-\ux$.

Then, we have $ q_{-}(t)\le q(t)\le q_+(t). $ Therefore, it
suffices to show that $q_{\pm}(t)$ achieves the value
$2D^{1/2-\lambda}$ for the first time at a moment
$\tilde\tau_{\pm}=\tilde s_n+o(1)$. Equation~\eqref{eqqpmeps}
yields
$$
\int\limits_{s_{n-1/2}+o(1)}^{\tilde\tau_{\pm}}\big(G_n(t)+1\pm\delta\big)dt=\int_{
2 Z_*(D)}^{\ox-\ux}\dfrac{dx}{F(x)}.
$$

Using the fact that $\tilde s_n$ is a root of~\eqref{eqsn} and
$s_{n-1/2}$ is a root of~\eqref{eqsn-1/2'} as well as the
definition of $\oF$ in~\eqref{eqoverF} and the symmetry of $F(x)$,
we have
$$
\begin{aligned}
\int\limits_{s_{n-1/2}}^{\tilde s_n}\big(G_n(t)+1\big)dt &
=\int\limits_{s_{n-1}}^{\tilde s_n}\big(G_n(t)+1\big)dt -
\int\limits_{s_{n-1}}^{s_{n-1/2}}\big(G_n(t)+1\big)dt\\
& =  \oF - \int_0^{\ox+\ux}\dfrac{dx}{F(x)} =
\int_0^{\ox-\ux}\dfrac{dx}{F(x)}.
\end{aligned}
$$
Hence,  the implicit function theorem in a neighborhood of
$(\tilde\tau_\pm,D)=(\tilde s_n,0)$  together with~\eqref{eqZ*o1}
yield $\tilde\tau_{\pm}=\tilde s_n+o(1)$.
\endproof

\subsection{Dynamics for $w\in (\ux,\ox-q_n)$ as $t\in(\tilde\tau_n,\tau_n)$}\label{subsubsec3}

 Assume that case~\ref{caseB} in the inductive step of the algorithm from section~\ref{subsecAlgorithm}
holds with some $s_n>\tilde s_n$. If case~\ref{caseA} holds, we
omit this step and proceed with section~\ref{subsubsec4} below.

The following lemma defines a time moment $\tau_n$
satisfying~\eqref{eqtimes2'}.

\begin{figure}[ht]
 \centering
\includegraphics[width=150mm]{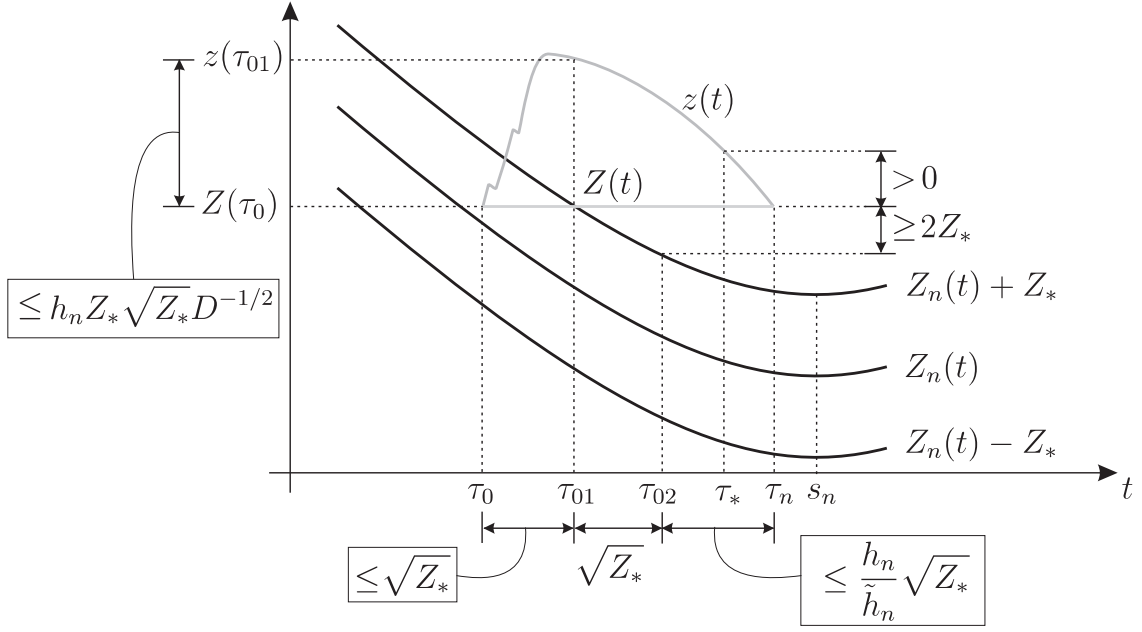}
\caption{Dynamics for $t$ close to $s_n$. Three black curves are
the graphs of $Z_n(t)$ and $Z_n(t)\pm Z_*$. The upper grey curve
is the graph  of the solution $z(t)$ to
equation~\eqref{eqdzdt2Phi}. The lower (horizontal)  grey curve is
the graph of $Z(t)$ given by~\eqref{eqZminz}.} \label{figznt}
\end{figure}

\begin{lemma}\label{lqGoesToSn}
There is a time moment $\tau_n=s_n+o(1)\le s_n$ such that $\dot
q(\tau_n)<0 $ and
\begin{equation}\label{eqqGoesToSn}
2(Z_n(\tau_n)-Z_*(D))D^{1/2}\le q(\tau_n)=Q(\tau_n)\le 2Z_*(D).
\end{equation}
\end{lemma}
\proof To follow the proof, we refer the reader to
Fig.~\ref{figznt}.
We consider equation~\eqref{eqdzdt2Phi} with the initial data
\begin{equation}\label{eqdzdt2PhiInitData}
 z(\tilde\tau_n)=Z_*(D)D^{-1/2}.
\end{equation}
 First,
note that $z(t)\ge Z(t)> Z_n(t)-Z_*(D)$ for all
$t\in[\tilde\tau_n,s_n]$. Indeed, if the equality
$Z(t)=Z_n(t)-Z_*(D) $ is achieved for the first time at some
moment $t\in [\tilde\tau_n,s_n]$, then $z(t)=Z(t)=Z_n(t)-Z_*(D)$.
Moreover, since $Z_n(t)-Z_*(D)$ is nonincreasing, we have $\dot
z(t)\le 0$ at this moment. On the other hand,
relations~\eqref{eqdzdt2Phi} and~\eqref{eqHnzeta} imply that $\dot
z(t)>0$. Thus, the first inequality in~\eqref{eqqGoesToSn} holds
with $\tau_n$ replaced by any $t\le s_n$.

Furthermore, $Z(t)$ is nonincreasing. Hence, $Z(t)\le
Z(\tilde\tau_n)\le Z_*(D)D^{-1/2}$ for all $t\ge\tilde\tau_n$ by
Lemma~\ref{ltau2s2}. It remains to find $\tau_n$ such that $z(\tau_n)=Z(\tau_n)$ and $\dot q(\tau_n)<0$.

Next, we introduce a positive function $a_*(D)$ such that
$a_*(D)\to 0$ and
\begin{equation}\label{eqz'NotSmall}
 \dot Z_n(t)\le - 2\sqrt{Z_*(D)},\quad t\le s_n-a_*(D).
\end{equation}
Consider the time moment
$$
\tau_*=\tau_*(D):=s_n-a_*(D)-\sqrt[4]{Z_*(D)}.
$$
First, assume
that
\begin{equation}\label{ddd}
z(\tau_*)>Z(\tau_*)
\end{equation}
(the case $z(\tau_*)=Z(\tau_*)$ will be considered at the end of the proof). Denote by $\tau_0$ the time
moment preceding $\tau_*$ such that
\begin{equation}\label{eqzbiggerZtau0*}
z(\tau_0)=Z(\tau_0),\qquad z(t)>Z(t)=Z(\tau_0),\quad
t\in[\tau_0,\tau_*].
\end{equation}
Obviously, $Z(\tau_0)\in[Z_n(\tau_0)-Z_*,Z_n(\tau_0)+Z_*]$.
Indeed, otherwise, $Z(\tau_0)>Z_n(\tau_0)+Z_*$ and, due
to~\eqref{eqHnzeta} and equation~\eqref{eqdzdt2Phi}, we would have
$\dot z(\tau_0)<0$.

Note that, due to~\eqref{eqz'NotSmall}, we have $Z(\tau_0)\ge
Z_n(\tau_0+\sqrt{Z_*})+Z_*$. Denote by
$$
\tau_{01}\in [\tau_0,\tau_0+\sqrt{Z_*}]
$$
a moment such that
\begin{equation}\label{eqZatTau01}
Z(\tau_0)=Z_n(\tau_{01})+Z_*.
\end{equation}
If $Z(\tau_{01})<Z(\tau_0)$, then there exists a desired time moment
$\tau_n\in[\tau_0,\tau_{01}]$ and the proof is complete.

Assume that $Z(\tau_{01})=Z(\tau_0)$. Then~\eqref{eqZatTau01}
yields
\begin{equation}\label{eqZatTau01'}
Z(\tau_{01})=Z_n(\tau_{01})+Z_*.
\end{equation}

 Let us estimate the velocity $\dot z(t)$ for
$t\in[\tau_0,\tau_{01}]$. For such $t$, we have $Z(t)\ge
Z_n(t)-Z_*$ and, due to~\eqref{eqZatTau01}, $Z(t)= Z(\tau_0)\le
Z_n(t)+Z_*$. Hence, \eqref{eqHnzeta} implies that $\dot z(t)\le
h_n Z_* D^{-1/2}$ for $t\in[\tau_0,\tau_{01}]$. Therefore,
\begin{equation}\label{eqztau01}
z(\tau_{01})\le Z(\tau_0)+h_n Z_* \sqrt{Z_*} D^{-1/2}.
\end{equation}

Set
$$
\tau_{02}:=\tau_{01}+\sqrt{Z_*}.
$$
From $\dot Z(t)<0$ it follows that $Z(t)>Z_n(t)+Z_*$ for $t\in[\tau_{01},\tau_{02}]$.
If $Z(\tau_{02})<Z(\tau_{01})$, then there exists a desired time
moment $\tau_n\in[\tau_{01},\tau_{02}]$ and the proof is complete.

Assume that $Z(\tau_{02})=Z(\tau_{01})$. It follows
from~\eqref{eqz'NotSmall} and~\eqref{eqZatTau01'} that
\begin{equation}\label{eqZtau02}
Z(\tau_{02})=Z(\tau_{01})=Z_n(\tau_{01})+Z_*\ge Z_n(\tau_{02})+Z_*
+ 2 Z_*.
\end{equation}
Furthermore, due to~\eqref{eqHnzeta} and
equation~\eqref{eqdzdt2Phi}, we have
$$
\dot z(t)<0,\quad t  \in[\tau_{01},\tau_{02}].
$$
Therefore (see~\eqref{eqztau01}),
\begin{equation}\label{eqztau02}
z(\tau_{02})\le z(\tau_{01})\le Z(\tau_0)+h_n Z_* \sqrt{Z_*}
D^{-1/2}.
\end{equation}

Now inequality~\eqref{eqZtau02}, relations~\eqref{eqHnzeta}, and
equation~\eqref{eqdzdt2Phi} imply that $ \dot z(t)\le -\tilde h_n
Z_* D^{-1/2} $ at least as long as $z(t)\ge
Z(\tau_{02})=Z(\tau_0)$, where $\tilde h_n$ does not depend on $t$
and $D$. Therefore, taking into account~\eqref{eqztau02}, we see
that $z(t)$ will achieve the value $Z(\tau_{02})=Z(\tau_0)$ at a
time moment $\tau_n\in [\tau_{02},\tau_{02}+h_n\sqrt{Z_*}/\tilde
h_n]$.
If follows from~\eqref{eqzbiggerZtau0*} and from the definitions
of $\tau_*$, $\tau_0$, $\tau_{01}$, $\tau_{02},$ and $\tau_n$ that
$\tau_n\in(\tau_*,s_n)$. Hence, $\tau_n$  is the desired time
moment.
%
%

Finally, if \eqref{ddd} is not valid, then $z(\tau_*)=Z(\tau_*)$
and, as we already know, $Z(\tau_*)\le  Z_* D^{-1/2}$. The
equality $z(\tau_*)=Z(\tau_*)$ implies $\dot q(\tau_*)\le 0$. If
$\dot q(\tau_*)<0$, then we can choose $\tau_n=\tau_*$. If $\dot
q(\tau_*)=0$, then relations \eqref{eqdzdt2Phi}, \eqref{eqHnzeta}
imply $Z(\tau_*)\in[Z_n(\tau_*)-Z_*,Z_n(\tau_*)+Z_*]$ and all the
argument following Eq.~\eqref{eqzbiggerZtau0*} can be repeated for
$\tau_0=\tau_*$.
\endproof

\subsection{Dynamics for $w\in (\ux,\ox-q_n)$ as $t\in(\tau_n,\theta_n)$}\label{subsubsec4}

%

\begin{lemma}\label{lsigma2tau2}
There is a time moment $ \theta_n=\tau_n+o(1)$ such that $\dot
q(t)<0$ for $t\in[\tau_n,\theta_n]$ and
$q(\theta_n)=2\big(y_n+o(1)\big)\big(Ds_n\big)^{1/2}.$
\end{lemma}
\proof 1. Setting $z(t)=\dfrac{q(t)}{2D^{1/2}}$, we consider the
problem (see Lemma~\ref{lqGoesToSn})
\begin{equation}\label{eqdzdt2Phi'}
\left\{
\begin{aligned}
& \dot z(t) = \dfrac{F(0)(H(t,Z(t))+o(1))}{2D^{1/2}},\\
& y_n s_n^{1/2}-Z_*(D)\le z(\tau_n)\le Z_* D^{-1/2},
\end{aligned}
\right.
\end{equation}
where $H(t,z)$ is given by~\eqref{eqHtz}. If $z(\tau_n)\le
Z_n(\tau_n)+Z_*(D)=y_n s_n^{1/2}+Z_*(D)+o(1)$, then we take
$\theta_n:=\tau_n$ and complete the proof.

Assume that $z(\tau_n)> Z_n(\tau_n)+Z_*(D)$.  We need to prove
that there exist positive functions $\varepsilon_*(D)$ and
$\delta_*(D)$ that tend to zero as $D\to 0$ and  such that the
solution $z(t)$ to~\eqref{eqdzdt2Phi} satisfies
\begin{equation*}
\dot z(t)<0,\quad t\in [\tau_n,\tau_n+\delta_*(D)],\qquad y_n\le
z(\tau_n+\delta_*(D))\le(y_n+\varepsilon_*(D))s_n^{1/2}.
\end{equation*}

To construct the functions  $\varepsilon_*(D)$ and $\delta_*(D)$,
we consider  positive sequences $\varepsilon_k,\delta_k\to 0$ such that $\delta_k\le s_n/2$. Due
to~\eqref{eq-Gn10}, \eqref{eqZn}, and~\eqref{eqFormulayn}, there
exists a sequence $c_k=c_k(\varepsilon_k)$ such that
\begin{equation}\label{eqdzdt2Phi''}
H(t,z)\le - c_k,\quad t\in\left[{\tilde s_n}, T\right],\ z\ge
Z_n(t)+\varepsilon_k.
\end{equation}
 Then~\eqref{eqdzdt2Phi'} and~\eqref{eqdzdt2Phi''} imply that
$z(t)$ achieves the value $Z_n(t)+\varepsilon_k$ at a time
moment $\theta_{nk}$ that satisfies
$$
0\le
\theta_{nk}-\tau_n\le\dfrac{Z_*(D)D^{-1/2}}{\inf\limits_{t\in[\tilde s,2s_n]}|\dot
z(t)|}\le \dfrac{4 Z_*(D)}{F(0)c_k}\le\delta_k
$$
for all $D\le D_k$, where $D_k$ is a strictly decreasing sequence
with $D_k\to0$. Therefore,
$z(\theta_{nk})=Z_n(\theta_{nk})+\varepsilon_k=(y_n+\tilde
\varepsilon_k)s_n^{1/2}$ where $\tilde \varepsilon_k\to 0$ as
$k\to\infty$. Now we set $\varepsilon_*(D):=\tilde\varepsilon_k$
and $\delta_*(D):=\delta_k$ for $D\in(D_{k},D_{k-1}]$.
\endproof

\subsection{Dynamics for $w\in (\ux,\ox-q_n)$ as $t\in(\theta_n,t_n)$}\label{subsubsec5}

 To describe the dynamics for
$t>\theta_n$, we consider~\eqref{eqdqdtnPhiMin} with the initial
data
\begin{equation}\label{eqInitialCondthetan}
q(\theta_n)=2\big(y_n+o(1)\big)\big(Ds_n\big)^{1/2}
\end{equation}
(see Lemma~\ref{lsigma2tau2}).
The following lemma will allow us to justify~\eqref{eqdotqrn1}
(with $n-1$ replaced by~$n$) and thus to complete the proof of
Theorem~\ref{tAsympImmortalFronts}.

\begin{lemma}\label{ltildethetan}
 There exists $\tilde\theta_n=s_n+o(1)\ge \theta_n$ and $r_n$ $(>s_n)$ not depending on $D$ such that the  solution
$q(t)$ to~$\eqref{eqdqdtnPhiMin}$, $\eqref{eqInitialCondthetan}$
satisfies $\dot q(t)>0$ for $t\in[\tilde\theta_n,r_n]$ and
$w(r_n)>-\ux$.
\end{lemma}
\proof Since $F(q)\ge 1/4-\ox^2$ for $q\in[0,\ox-\ux]$, it
suffices to estimate the expression in brackets
in~\eqref{eqdqdtnPhiMin}. Using the definition of $Q(t)$, the
monotonicity of $E(\cdot)$, and
relations~\eqref{eqInitialCondthetan} and~\eqref{eqGnsn1} (with
$n-1$ replaced by $n$), we see that there exists
$\tilde\theta_n=s_n+o(1)\ge \theta_n$ such that
\begin{equation}\label{eqdqdtnPhiMin1}
\begin{aligned}
& -G_n(t)- 2E\left(\dfrac{Q(t)}{2(D
t)^{1/2}}\right)+1+\mu(t,D)\\
&\qquad \ge -G_n(t)- 2E\left(\dfrac{q(\theta_n)}{2(D
t)^{1/2}}\right)+1+\mu(t,D)\\
&\qquad =G_{n+1}(t)+1+o(1)>0\quad\text{for }  t\ge \tilde\theta_n.
\end{aligned}
\end{equation}

Equation~\eqref{eqdqdtnPhiMin} and
inequality~\eqref{eqdqdtnPhiMin1} imply that $\dot q(t)>0$ for
$t\ge \tilde\theta_n.$ Obviously, we can choose $r_n>s_n$
independent of $D$ such that $w(r_n)=\ox-q(r_n)>-\ux$.
\endproof

Finally, we determine the time moment $t_n$ (see~\eqref{eqtimes3})
at which the $n$-th front becomes steady and immortal and its
position $x_n=\ox-q_n$ as follows:
$$
t_n:=\sup\{t\in [\theta_n,\tilde\theta_n]: \dot q(t)=0,\
q(t)=Q(t)\},\qquad q_n:=q(t_n)=Q(t_n).
$$
Due to Lemmas~\ref{lsigma2tau2} and~\ref{ltildethetan},
$t_n=s_n+o(1)$ and
\begin{equation}\label{eqdotqn0tn}
\dot q(t)>0\quad\text{for } t\in(t_n,r_n].
\end{equation}
 This
justifies~\eqref{eqdotqrn1} (with $n-1$ replaced by $n$).
Moreover, since $\dot q(t_n)=0$, the value $y:=\dfrac{Q(t_n)}{2(D
t_n)^{1/2}}$ satisfies (due to~\eqref{eqdqdtnPhiMin}
and~\eqref{eqMuoD12})
$$
0=-G_n(s_n+o(1))- 2E(y)+1+o(1).
$$
On the other hand,   the definition~\eqref{eqHtz} of $H_n(t,z)$,
the definition~\eqref{eqZn} of $Z_n(t)$ and the
definition~\eqref{eqFormulayn} of $y_n$ imply that $y_n$ is the
root of the equation
$$
0=-G_n(s_n)- 2E(y_n)+1.
$$
 Therefore,  $y=y_n+o(1)$, which implies
$q(t_n)=Q(t_n)=2\big(y_n+o(1)\big)\big(Ds_n\big)^{1/2}$. This
completes the proof of Theorem~\ref{tAsympImmortalFronts}.

\begin{remark}\label{remn1}
An essential ingredient in the proof is the
approximation~\eqref{eqUcloseE}  of the integral of the solution
$U$ by the error function $E(\cdot)$. It is valid for $t$
separated from zero and thus can be used for $n\ge 2$. In the case
$n=1$, the analog of Eq.~\eqref{eqdqdtn-1/2Phi} will be $ \dot
q(t)=f(q)\uU(t)$ on the interval $t\in(0,t_{1/2})$, where
$q(t)=\ox-w(t)$, $q(0)=o(1)$, and $q(t_{1/2})=\ox+\ux$. On this
time interval, one can use the approximation $\uU(t)=1+o(1)$
(see~\eqref{eqIC} and~\eqref{limit''}) instead
of~\eqref{eqUcloseE}.
\end{remark}

\section{Numerics}

\begin{figure}[ht]
 \centering
\includegraphics[width=160mm]{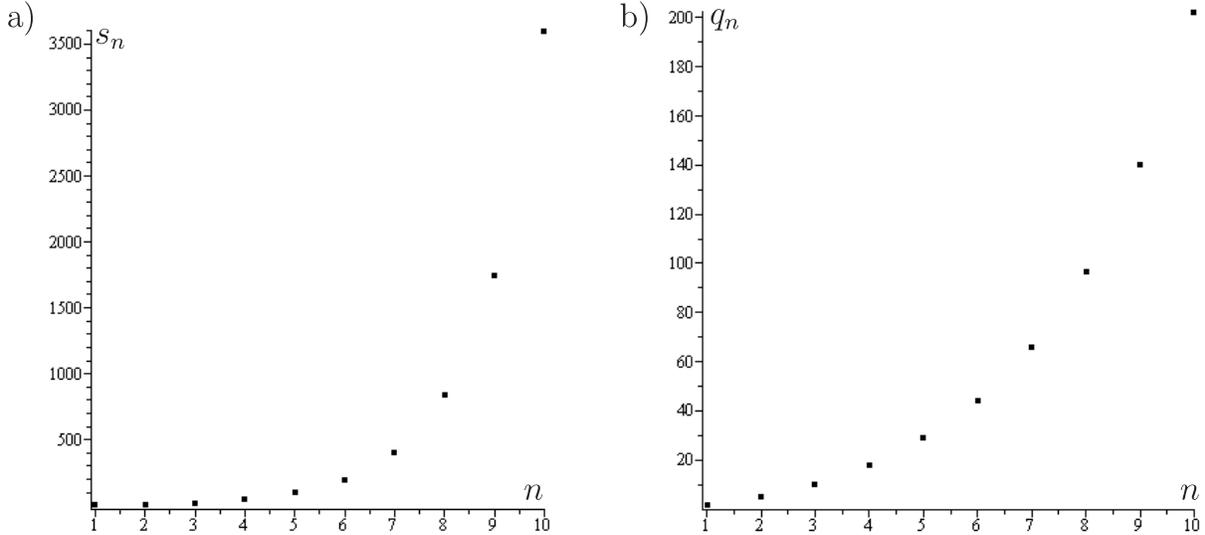}
 \caption{The
two graphs indicate the values of a) $s_n$ and b) $q_n$ for
$n=1,\dots,10$ from Theorem~\ref{tAsympImmortalFronts}. The
maximal and minimal values for admissible thresholds are
$\ux=1/100$ and $\ox=1/4$, respectively.} \label{figGraph}
\end{figure}

Figure~\ref{figGraph} illustrates the  values of $s_n$ and $q_n$
for $n=1,\dots,10$ from Theorem~\ref{tAsympImmortalFronts} found
numerically.
We   note that $\tilde s_n=s_n$ for $n=1,\dots,6$ (i.e.,
case~\ref{caseA} from the inductive step of the algorithm from section~\ref{subsecAlgorithm} takes place) and $\tilde
s_n < s_n$ for $n=7,\dots,10$ (i.e., case~\ref{caseB} from the
inductive step takes place). In general, it is an open question
whether one of the two  cases~\ref{caseA} and~\ref{caseB} takes place for all $n$, or both
may fail for some $n$ and hence the resulting
sequences~\eqref{eqSequencesnyn} are finite.

\section*{Acknowledgments}
Pavel Gurevich acknowledges the support of the DFG through the Collaborative Research
Center 910 and the Heisenberg fellowship. Dmitrii Rachinskii was supported by National
Science Foundation, grant DMS-1413223. The authors are grateful to
Sergey Tikhomirov who created a software for a number of numerical
experiments.

\end{document}